  \documentclass[]{amsart}
   \newtheorem{lemma}{Lemma}[section]
   \newtheorem{theorem}[lemma]{Theorem}

   \newtheorem{coro}[lemma]{Corollary}
   \newtheorem{definition}[lemma]{Definition}
   \newcommand{\eps}{\varepsilon}
   \newcommand{\Wsob}{\smash{{\stackrel{\circ}{W}}}_2^1(D)}

\newcommand{\D}{\Delta}
\newcommand{\p}{\partial}
\renewcommand{\phi}{\varphi}

\renewcommand{\k}{\kappa}

\title[Stochastic Coupled Atmosphere--Ocean Model]
{Stochastic Dynamics\\ of a Coupled Atmosphere--Ocean Model}

\author{Jinqiao Duan}
\address[J. Duan]
{Department of Applied Mathematics\\
 Illinois Institute of Technology\\
   Chicago, IL 60616, USA }
\email[J.~Duan]{duan@iit.edu}

\author{Hongjun Gao}
\address[H. Gao]
{Department of Mathematics\\
Nanjing Normal University\\
Nanjing 210097, China}
\email[H.~Gao]{gaohj@njnu.edu.cn}

\author{Bj{\"o}rn Schmalfu{\ss}}
\address[B. Schmalfu{\ss}]{%
  Department of Sciences\\
  University of Applied Sciences\\
  Geusaer Stra{\ss}e\\
  06217 Merseburg, Germany}
\email[B.~Schmalfu{\ss}]{bjoern.schmalfuss@in.fh-merseburg.de}

\date{February 5, 2002}

\subjclass{Primary 60H15;  Secondary  86A05, 34D35}
\keywords{Stochastic  geophysical flow models, random attractor, climate dynamics,  finite degrees of freedom,  ergodicity}

\begin{document}

\begin{abstract}

{\bf Stochastics and Dynamics}. Vol 2, No. 3, 2002,   357-380.

The investigation of the   coupled atmosphere-ocean system
is not only scientifically challenging but also practically important.

We consider a coupled atmosphere-ocean model, which involves hydrodynamics,
thermodynamics, and random  atmospheric dynamics due to short time
influences
  at the air-sea interface.  We reformulate this model as a random
dynamical system.
First,  we have shown that the asymptotic dynamics
of the coupled atmosphere-ocean model  is described by a
random climatic attractor.
Second, we have estimated the
atmospheric temperature  evolution  under oceanic feedback,
 in terms of the freshwater flux, heat flux
and the   external fluctuation at the air-sea interface, as well as
the earth's longwave radiation coefficient
and the shortwave solar radiation profile.
Third, we have demonstrated that this system has finite
degree of freedom by presenting a finite set of determining
functionals in probability.
Finally, we have proved that the coupled atmosphere-ocean model
is ergodic under suitable conditions for physical parameters and randomness,
and thus for any observable of
the    coupled  atmosphere-ocean flows, its time average
approximates the statistical ensemble average, as long as  the time
interval is sufficiently long.

\end{abstract}

\maketitle

\section{Geophysical background}

The coupled atmosphere-ocean system defines the environment we live.
Randomness or uncertainty is ubiquitous in this coupled, complex, multiscale system: For example,
stochastic forcing (wind stress, heat flux and freshwater flux),
uncertain parameters, random sources or inputs, and
random  boundary conditions.

Mathematical models are a key component of our
understanding of  climate and geophysical  systems.
It is our belief that the fidelity of these
models to nature can greatly benefit through the inclusion of stochastic
 effects.
Taking stochastic effects into account is of central
   importance for the  development of mathematical models of many
   phenomena in geophysical  and climate
flows.

We   consider a two-dimensional
coupled  atmosphere-ocean model in the
latitude-depth  plane, with   atmospheric dynamics highly
simplified, i.e., the atmospheric dynamics is described by an energy
 balance model.  The oceanic dynamics is described by the Navier-Stokes
equation in vorticity form  and the transport equations for  heat and salinity.
The energy balance model is under random impact due to,
 for example, eddy transport fluctuation, stormy
bursts of latent heat, and flickering
 cloudiness variables.
So this coupled  atmosphere-ocean
model consists of stochastic and
deterministic partial differential equations, together with air-sea flux or
Neumann boundary conditions.
We will reformulate this model as a random dynamical system.

\medskip

The ocean and the atmosphere are constantly interacting through
the air-sea exchange process. The ocean moves much slower than
the  atmosphere does.
It is generally believed that the   ocean
plays an important role in the global climate dynamics in relatively long time scales, due to ocean's large capacity of holding and transporting huge amount of heat or cold
around the globe \cite{SieChuGou01} . However, a complete quantitative understanding or estimate
for ocean's impact on climate is  lacking.  A particular issue is: How does
the ocean affect or provide feedback to the air temperature, which is the most
important climate quantity we are usually concerned about?
 This is a direct impact of
 the ocean on the climate. It is desirable to predict
or estimate this feedback in the context of our simple coupled
atmosphere-ocean model.

The existence and interpretation of climatic attractors have been
controversial and have caused a lot of debate \cite{NicNic84}. A
low dimensional climatic attractor was regarded as an indication
that the main feature of long-time climatic evolution may be
viewed as the manifestation of a {\em deterministic} dynamics.
We will consider {\em random} climate attractors,  and   the long
time regimes that such attractors   represent still carry the
stochastic information of the geophysical flow system.
We will also investigate the finite dimensionality of the asymptotic
dynamics by checking the {\em determining functionals in probability}.

In a special case of physical parameters and random noise,
we obtain a random attractor which is defined by a single
random variable. This  random variable attracts  all other
motions exponentially fast. This random variable corresponds to a unique
invariant measure, which is the  expectation of the
Dirac measure with the random variable
as the random mass point; see \cite{Arn98}.  In this case, the coupled
atmosphere-ocean model is {\em ergodic}, and thus for any observable of
the    coupled  atmosphere-ocean flows, its time average
approximates the statistical ensemble average, as long as  the time
interval is sufficiently long.

\bigskip

In the  next section, we present the coupled
atmosphere-ocean model,  and discuss the well-posedness
of this coupled   model in \S 3. Then we investigate
the dissipativity property   in \S 4.  This   property is the basis of
the asymptotic behavior  of the coupled system to be considered in \S 5:
  atmospheric temperature  evolution (with oceanic feedback),
random attractors, finite dimensionality and ergodicity.
 Finally, we summarize these results  in   \S 6.

\section{A coupled atmosphere--ocean model}\label{sect2}

We   consider a zonally averaged,
coupled atmosphere-ocean model on the meridional,
latitude-depth $(y,z)$-plane as used by various authors
\cite{StoWriMys92, WriSto91, CheGhi96, DijNee00, Dij00}.
It is composed of a one-dimensional
 stochastic energy
balance model   proposed by North and Cahalan \cite{NorCah81},
for the    latitudinal
 atmosphere surface temperature  $\Theta(y,t)$ on domain $0<y <1$, together with
the Boussinesq equations for ocean dynamics in terms of vorticity $q(y,z,t)$,
and transport  equations for the oceanic salinity $S(y, z, t)$
and the oceanic temperature $T(y, z, t)$ on
the domain $D=\{ (y,z): 0\leq y, z \leq 1\}$:
\begin{align}\label{eqn1}
\begin{split}
\Theta_t   = &  \Theta_{yy} -(a +\Theta) +
S_a(y)-b(y)(S_o(y)+\Theta -T(y,1))
        + \dot{w},  \\
q_t + J(q, \psi )      = & {\rm Pr} \D q  +{\rm Pr  \cdot Ra} (\partial_yT - \partial_yS_y), \\
T_t + J(T, \psi )        = & \D T, \\
S_t + J(S, \psi )          = & \D S,
\end{split}
\end{align}
where
\[
q(y, z, t)= -\D \psi
\]
is the vorticity, $a$ is a   positive constant
parameterizing the effect
of the earth's longwave radiative cooling, $S_a(y)$ and $ S_o(y)$
are empirical functions representing the latitudinal
dependence of the shortwave solar radiation,
$b(y)$  is  the latitudinal  fraction of the earth covered
by the ocean basin,  Pr is the
Prandtl number  and Ra is the Rayleigh number.  The first equation is the
energy balance model proposed by
North and Cahalan \cite{NorCah81}. The  fluctuating  forcing
$\dot{w}(y, t)$ may arise from, for example, eddy transport fluctuation, stormy
bursts of latent heat, and flickering
 cloudiness variables.  This forcing term is
usually  of a shorter time scale than the response time scale of
the large scale oceanic thermohaline circulation. So we neglect the
autocorrelation time of this fluctuating process as in \cite{NorCah81}.
We thus assume the noise is white in time.
The  random white-in-time noise  $\dot{w}(y, t)$
is described as the generalized time derivative of a Wiener process $w(y, t)$
with mean zero and covariance operator $Q$.
Moreover, $J(g,h)=g_xh_y-g_yh_x$
is the Jacobian operator and $\D=\p_{yy}+\p_{zz}$ is the Laplacian
operator.   All these equations are in non-dimensionalized forms.

Note that the Laplacian operator $\D$ in the temperature and salinity transport
equations  is presumably $\p_{yy}+ \frac{\k_H}{\k_V} \delta^2 \p_{zz}$
with  $\delta$ being the aspect ratio,  and $\k_H, \k_V$ the horizontal
and vertical diffusivities of heat/salt, respectively.  However, our
energy-type estimates and the results below will not be essentially
affected by taking a homogenized  Laplacian operator $\D=\p_{yy}+\p_{zz}$.
All  our results would be true for this modified Laplacian.
The effect of the rotation   is
parameterized in the magnitude of the viscosity and diffusivity terms
as discussed in \cite{ThuMcW92}.

The  no-flux boundary condition  is taken for
the  atmosphere  temperature $\Theta (y, t)$
\begin{equation*}
\Theta_y (0, t)=\Theta_y (1, t)  = 0.
\end{equation*}
The  fluid boundary condition is no normal flow and free-slip
on the whole boundary
\begin{equation*}
\psi=0,  \;  q = 0.
\end{equation*}

The flux boundary conditions are assumed
for the ocean temperature $T$ and salinity $S$.

 At top $z=1$, the fluxes are specified as:
\begin{equation}
\label{btop}
 \partial_z T(y,1) = S_o(y)+ \Theta(y)-T(y,1), \; \;    \partial_z S(y,1)=F(y),
 \end{equation}
with $F(y)$ being the given freshwater flux.

At bottom $z=0$:
\begin{equation*}
\partial_z T =  \partial_z S=0  \;.
\end{equation*}

On the lateral boundary    $y\in\{0, 1\}$:
\begin{equation*}
\label{side} \partial_y T=\partial_y S =0.
\end{equation*}


The stochastic partial differential equation for air temperature $\Theta$
in (\ref{eqn1})
is only defined on the air-sea
interface ($0\leq y \leq 1$) and it may be regarded as a   dynamical
boundary condition.  The boundary condition (\ref{btop}) involves
a coupling between the atmospheric and oceanic temperature at the
air-sea interface.

The deterministic version of this model was studied in \cite{gaodua02}.
Now we look at the well-posedness of this coupled atmosphere-ocean model
and then investigate its random dynamics.

\section{Well-posedness }

In this section we will show that (\ref{eqn1}) defines a well posed
model. In particular, we can show that (\ref{eqn1}) has a unique solution.
Without such a  property it would not be possible to make predictions
 from the model
numerical simulations or investigate the stability behavior.\\
Now we are going to re-formulate the model such that appropriate
tools of the theory of random dynamical system can be applied to analyse
the coupled atmosphere-ocean model under a random wind forcing.
For the following we need some tools from the theory of partial differential equations.\\
Let $W_2^1(D)$ be the Sobolev space of functions on $D$ with first
generalized derivative in $L_2(D)$, the function space of square
integrable functions on $D$ with norm and  inner product
\[
\|u\|_{L_2}=\left(\int_D|u(x)|^2dD\right)^\frac{1}{2},\quad
(u,v)_{L_2}=\int_Du(x)v(x)dD,\quad u,\,v\in L_2(D).
\]
The space $W_2^1(D)$ is equipped with the norm
\[
\|u\|_{W_2^1}=\|u\|_{L_2}+\|\partial_y u\|_{L_2}+\|\partial_z u\|_{L_2}.
\]
Motivated by the zero-boundary conditions of $q$ we also introduce the space
$\Wsob$ which contains roughly speaking  functions which are zero on the boundary $\partial D$
of $D$.
This space can be equipped with the norm
\begin{equation}\label{eqno}
\|u\|_{\Wsob}=\|\partial_y u\|_{L_2}+\|\partial_z u\|_{L_2}.
\end{equation}

Similarly, we can define function spaces on the interval $(0,1)$
denoted by $L_2(0,1)$ and $W_2^1(0, 1)$.\\
Another Sobolev space is given by $\dot{W}_1^2(D)$ which is a
subspace of $W_2^1(D)$ consisting of functions $u$ such that
$\int_DudD=0$. A norm equivalent to the $W_2^1$-norm on
$\dot{W}_1^2(D)$ is given by the right hand side of (\ref{eqno}).
For functions in $L_2(D)$
having this property we will write $\dot{L}_2(D)$.\\

Since $D$ has a Lipschitz continuous boundary $\partial D$ there exists a {\em continuous
trace operator}:
\begin{equation*}
\gamma_{\partial D}: W_2^1(D)\to H^\frac{1}{2}(\partial D).
\end{equation*}
Here $H^\frac{1}{2}(\partial D)$ is a boundary space, see Adams \cite{Ada75}
or below.
Similarly, we can introduce trace operators that map onto a part of the boundary of $\partial D$
for instance for the subset $\{(y,z)\in \bar D|z=1\}$ of $\bar D$.
For this mapping we will write
\begin{equation}\label{trace1}
\gamma_{z=1}: W_2^1(D)\to H^\frac{1}{2}(0,1).
\end{equation}
The adjoint operator
\[
\gamma_{z=1}^\ast:(H^\frac{1}{2}(0,1))^\prime\to ({W}_2^1(D))^\prime
\]
is  also continuous. Note that
$^\prime$ denotes the dual space for a given Banach space.\\

Our intention is now to formulate the problem (\ref{eqn1}) with the non--homogeneous boundary
conditions in a weak form.
For convenience,  we introduce the vector notation for unknown
geophysical quantities
\begin{equation}\label{eqquant}
u=(\Theta,q,T,S).
\end{equation}

\bigskip

We now take the linear differential operator from (\ref{eqn1}):
\begin{equation*}
{\mathcal{A}}u =\left(
\begin{array}{l}
-\partial_{yy}^2\Theta+(1 + b(y)) \Theta\\
-{\rm Pr} \Delta q\\
-\Delta T\\
-\Delta S
\end{array}
\right).
\end{equation*}
Remind that the function $1\ge b(y)\ge 0$. ${\mathcal{A}}$ is
defined on functions that are sufficiently smooth. We also have
the following boundary conditions from Section \ref{sect2}
\begin{align*}
\begin{split}
&\partial_y\Theta(0)=\partial_y\Theta(1)=0,\\
&\psi|_{\partial D}=0,\quad q|_{\partial D}=0,\\
& \partial_zT(y,1)=S_o(y)+\Theta(y)-T(y,1),\quad \partial_zS(y,1)=F(y),\\
&\partial_z T(y,0)=\partial_z S(y,0)=0,\\
&\partial_y T(0,z)=\partial_y S(0,z)=0,\quad \partial_y T(1,z)=\partial_z S(1,z)=0.
\end{split}
\end{align*}

We will assume that $S_o,\,S_a$ and $F\in L_2(0,1)$.
Note that
\[\frac{d}{dt}\int_{\Omega} S dydz = \int_0^1 F(y) dy = \mbox{constant}.
\]
It is reasonable (see \cite{Dij00})  to assume that
\begin{equation}
\int_0^1 F(y) dy = 0,
\end{equation}
and thus $ \int_DSdydz$ is constant in time
and we may assume it is zero:
\begin{equation*}
 \int_DSdydz = 0.
\end{equation*}
Thus we have the usual Poincar\'e  inequality for $S$.
Unfortunately, this is not the situation for $T$. However,
 we can derive the following Poincar\'e  inequality
\begin{equation}
\label{tcoest} \|T\|^2 \le  2\|\gamma_{z=1}T\|_{L^2}^2 + 4\|\nabla T\|^2,
\end{equation}
as in Temam \cite{Tem97}, Page 52.

We  introduce the phase space for our geophysical quantities
$H=L_2(0,1)\times L_2(D)\times L_2(D)\times \dot{L}_2(D)$ with the usual $L_2$ inner product and
$V=W_2^1(0,1)\times \Wsob\times W_2^1(D)\times \dot{W}_2^1(D)$.
For another sufficiently smooth functions
$v=(\bar \Theta,\bar q,\bar T,\bar S)$,
we can calculate via integration by parts
\begin{align}\label{eq1}
\begin{split}
({\mathcal{A}}u,v)_{H}
=&\int_{(0,1)}\partial_{y}\Theta\partial_{y}\bar \Theta dy
+
\int_0^1(1 + b)\Theta\bar\Theta dy\\
&+
{\rm Pr}\int_D\nabla q\cdot\nabla \bar q\,dD \\
&
+
\int_D\nabla T\cdot\nabla \bar T\,dD-\int_0^1(S_o(y)+
\Theta(y)-\gamma_{z=1}T(y,z))(\gamma_{z=1}\bar T(y,z))dy\\
&
+
\int_D\nabla S\cdot\nabla \bar S\,dD-\int_0^1F(y)(\gamma_{z=1}\bar S(y,z))dy.
\end{split}
\end{align}
Hence on the  space $V$ we can introduce
a bilinear form $\tilde{a}(\cdot, \cdot)$
which is continuous, symmetric and positive
\begin{align*}
\tilde{a}(u,v)=&\int_{(0,1)}\partial_{y}\Theta\cdot\partial_{y}\bar \Theta dy
+ \int_0^1(1 + b)\Theta\bar\Theta dy
+
{\rm Pr}\int_D\nabla q\cdot\nabla \bar qdD \\
&
\\
&
+
\int_D\nabla T\cdot\nabla \bar TdD
+
c_0\int_0^1\gamma_{z=1}T \gamma_{z=1}\bar Tdy
+
\int_D\nabla S\cdot\nabla \bar SdD
\end{align*}
for some $c_0>0$.
The other terms from (\ref{eq1}) will be considered separately.
This bilinear form defines a unique linear continuous operator $A:V\to V^\prime$ such that
$\langle Au,v\rangle=\tilde{a}(u,v)$.
We can  see that   the bilinear form  $\tilde{a}(u,v)$ is positive,    using the Poincar\'e inequality
(\ref{tcoest}).
According to (\ref{eq1}), we  now introduce the nonlinear operator $F(u):=F_1(u)+F_2(u)$ where
\[
F_1(u)[y,z]=
\left(
\begin{array}{c}
0\\
-J(q,\psi)\\
-J(T,\psi)\\
-J(S,\psi)
\end{array}
\right)[y,z].
\]
and
\[
F_2(u)[y,z]=
\left(
\begin{array}{l}
-a+S_a(y)-b(y)(S_o(y)-\gamma_{z=1}T)\\
{\rm Pr}\,{\rm Ra}(\partial_yT-\partial_y S)\\
\gamma_{z=1}^\ast(S_o(y)+\Theta-\gamma_{z=1}T) \\
\gamma_{z=1}^\ast F(y)
\end{array}
\right)[y,z].
\]

\begin{lemma}\label{lcon}
The operator $F_1:V\to H$ is continuous. In particular, we have
\[
\langle F_1(u),u\rangle=0.
\]
\end{lemma}
\begin{proof}
We have a constant $c_1>0$ such that
\begin{equation}\label{eq6a}
\|\psi\|_{W_2^3(D)}\le c_1\|q\|_{W_2^1(D)}
\end{equation}
for any $q\in W_2^1(D)$ which follows straightforwardly by regularity properties of
a linear elliptic boundary problem. Note that
$W_2^3$ is a Sobolev space with respect to the third derivatives.
Hence we get:
\begin{equation*}
\begin{split}
\|J(T,\psi)\|_{L_2} &
\le \sup_{(y,z)\in D}(|\partial_y\psi(y,z)|+|\partial_z\psi(y,z)|)\times\\
&
\times\left(\int_D|\partial_yT(y,z)|+|\partial_zT(y,z)|dD\right).
\end{split}
\end{equation*}
  The second factor on the right hand side is bounded by
\[
\left(\int_D|\partial_yT(y,z)|^2dD\right)^\frac{1}{2}
+\left(\int_D|\partial_zT(y,z)|^2dD\right)^\frac{1}{2}\le \|u\|_V.
\]
On account of the Sobolev embedding
Lemma,  we have some  positive constants  $c_2,\, c_3$ such that
\[
\sup_{(y,z)\in D}(|\partial_y\psi(y,z)|+|\partial_z\psi(y,z)|)
\le c_2 \|\nabla \psi\|_{W_2^2(D)}\le c_3\|q\|_{W_2^1(D)}\le c_3 \|u\|_V.
\]
Hence we have a positive constant $c_4$ such that
\[
\|J(T,\psi)\|_{L_2}\le c_4\|u\|_V^2
\]
for $u\in V$.
Similarly, we can treat the other terms containing $J$.\\
We now show that
\[
\langle J(T,\psi),T\rangle =0.
\]
For the other terms containing $J$ we get a similar property. We obtain via
  integration by parts
\begin{align*}
\int_D\partial_yT&\partial_z\psi\, TdD-\int_D\partial_zT\partial_y\psi\, TdD\\
=&-\int_D\partial_{yz}^2T\psi \,TdD+\int_D\partial_{zy}^2T\psi T\,dD
-\int_D\partial_y T\psi\partial_zT\,dD+\int_D\partial_z T\psi\partial_yT\,dD\\
&+\int_{(0,1)}\partial_yT\psi T|_{z=0}^{z=1}dy-
\int_{(0,1)}\partial_zT\psi T|_{y=0}^{y=1}dz=0
\end{align*}
because $\psi$ is zero on the boundary $\partial D$. This relation
is true for a set of  sufficiently smooth functions $\psi,\,T$
which are dense in $\Wsob\times W_2^1(D)$. By the continuity of
$F_1$, as just shown in Lemma \ref{lcon}, we can extend this
property to $\Wsob\times W_2^1(D)$.
\end{proof}
\begin{lemma}
The following estimate holds
\[
\|F_2(u)\|_{V^\prime}\le c_5\|u\|_V+c_6.
\]
 for some  positive constants $c_5,\,c_6$.
\end{lemma}
\begin{proof}
Let $\zeta\in W_2^1(D)$.
Since $\gamma_{z=1}T\in H^\frac{1}{2}(0,1)\subset H^{-\frac{1}{2}}(0,1)$,
we have $\gamma_{z=1}^\ast\gamma_{z=1}T\in (W_1^2(D))^\prime$ and
\[
|\langle \gamma_{z=1}^\ast\gamma_{z=1}T,\zeta\rangle|
=
|\langle \gamma_{z=1}T,\gamma_{z=1}\zeta\rangle|\le c_7\|T\|_{W_2^1}\|\zeta\|_{W_2^1}
\]
for $c_7>0$,
which immediately gives the first part of the above inequality. The other parts can be
treated similarly.
\end{proof}

After this preparation, we are able to write our problem as a stochastic evolution equation.
An introduction into the theory of stochastic differential equations one can be found
in Zabczyk \cite{Zab00}.
\\
Let $\dot{w}$ be a noise on $L_2(D)$ with finite energy  given by the
covariance operator  $Q$ of the Wiener process $w(t)$ which is defined on a probability space
$(\Omega,\mathcal{F},\mathbb{P})$.
For the vector
\[
W=(w,0,0,0)
\]
we rewrite the coupled atmosphere-ocean system (\ref{eqn1})
as a     stochastic differential equation on   $V^\prime$:
\begin{equation}\label{eq19}
\frac{du}{dt}+Au=F(u)+\dot{W},\qquad u(0)=u_0\in H,
\end{equation}
where $\dot{w}$ is a white noise as the generalized temporal derivative
of a Wiener process  $w$  with continuous trajectories on $\mathbb{R}$  and with
values in $L_2(0,1)$.
Sufficient for this regularity
is that the trace of the covariance is finite with respect to the space $L_2(0,1)$:
${\rm tr}_{L_2}Q<\infty$.
In particular, we can choose the canonical probability space where the set of elementary
events $\Omega$ consists of the paths of $w$ and the probability measure $\mathbb{P}$
is the Wiener measure with respect to covariance $Q$.
\\

In the following, we need  a stationary Ornstein-Uhlenbeck process solving the
linear  stochastic equation on $(0,1)$
\begin{equation}\label{eq20}
\frac{dz}{dt} + A_1z = \dot{w}
\end{equation}
where $A_1 = - \partial_{yy} + ( 1 + b(y))$ is the linear operator  with the homogeneous Neumann
boundary condition at $y =0 $ and $y=1$.

\begin{lemma}
Suppose that the covariance $Q$ has a finite trace : ${\rm tr}_{L_2} Q<\infty$. Then
(\ref{eq20}) has a unique stationary solution generated by
\[
(t,\omega)\to z(\theta_t\omega).
\]
Moreover,  $Z(\omega)=(z(\omega),0,0,0)$ is a random variable in $V$.
\end{lemma}

For the proof we refer to Da Prato and Zabczyk \cite{DaPZab92}, Chapter 5,
or Chueshov and Scheutzow \cite{ChuScheu01}.
\\

For our calculations it will be appropriate to transform (\ref{eq19})
into a differential equation without white noise but with random coefficients.
We  set
\begin{equation}\label{eq19a}
v:=u-Z
\end{equation}
Thus we obtain a random differential equation in $V^\prime$
\begin{equation}\label{eq18}
\frac{dv}{dt}+Av=F_1(v)+F_2(v+Z(\theta_t\omega)),\qquad v(0)=v_0\in H .
\end{equation}
Equivalently, we can formulate the  equation (\ref{eq18}) using test functions
\begin{equation*}
\frac{d}{dt}(v(t), \zeta) + a(v(t), \zeta) = (F_1(v(t)), \zeta)
+ (F_2(v(t)+Z(\theta_t\omega)), \zeta)\quad \text{for all } \zeta\in V.
\end{equation*}

We have obtained a differential equation without white noise but
with random coefficients.
Such a differential equation can be treated sample-wise
for {\em any} sample $\omega$. Hence it
  is simpler to consider (\ref{eq18})
  than to study the stochastic differential equation (\ref{eq19})
  directly.
We are looking for solutions in
\begin{equation*}
v\in C([0, \tau]; H)\cap L^2(0, \tau; V),
\end{equation*}
for all $\tau>0$. If we can solve this equation then $u:=v+Z$
defines a solution version of (\ref{eq19}). For the well posedness
of the problem we now have the following result.

\begin{theorem}\label{tEX}
({\bf Well-Posedness})
For any time $\tau>0$, there exists a unique solution of (\ref{eq18})
in $C([0,\tau];H)\cap L_2(0,\tau;V)$. In particular, the solution mapping
\[
{{\mathbb{R}}}^+\times \Omega\times H\ni(t,\omega,v_0)\to v(t)\in H
\]

is measurable in its arguments and the solution mapping
$H\ni v_0\to v(t)\in H$ is continuous.\\
\end{theorem}
\begin{proof}
By the properties of $A$ and $F_1$ (see Lemma \ref{lcon}),
the random   differential equation  (\ref{eq18})  is   essentially similar  to the
2 dimensional  Navier Stokes equation.  Note that $F_2$ is only an affine mapping.
Hence we have existence and uniqueness
and the above regularity assertions.
\end{proof}

On account of the transformation (\ref{eq19a}),  we find that (\ref{eq19}) also has a unique solution.

Since the solution mapping
\[
{\mathbb{R}}^+\times\Omega\times H\ni (t,\omega,v_0)\to v(t,\omega,v_0)=:\phi(t,\omega,v_0)\in H
\]
is well defined,  we can introduce a random dynamical system. On $\Omega$ we can define
a shift operator $\theta_t$ on the paths of the Wiener process
that pushes our noise:
\[
w(\cdot,\theta_t\omega)=w(\cdot+t,\omega)-w(t,\omega)\quad \text{for }t\in{\mathbb{R}}
\]
which is called the {\em Wiener shift}.
Then $\{\theta_t\}_{t\in{\mathbb{R}}}$ forms a flow which is ergodic for the probability measure
${\mathbb{P}}$.
The properties of the solution mapping cause  the following relations
\begin{align*}
&\phi(t+\tau,\omega,u)=\phi(t,\theta_\tau\omega,\phi(\tau,\omega,u))\quad\text{for }t,\,\tau\ge 0\\
&\phi(0,\omega,u)=u
\end{align*}
for any $\omega\in\Omega$ and $u\in H$. This property is called the cocycle property of $\phi$
which is important to study the dynamics of random systems. It is a generalization
of the semigroup property. The cocycle $\phi$ together with the flow $\theta$ forms
a {\em random dynamical system}.

\section{Dissipativity}

In this section we are going to show that   the coupled  atmosphere-ocean system (\ref{eqn1})
 is dissipative,  in the sense that it has an absorbing  (random) set.
 This definition has been used for deterministic systems
\cite{Tem97}.
This means that the solution vector
$v$ is contained in a particular region of the phase space $H$ after a sufficiently long time.
Dissipativity will be very important for understanding  the
asymptotic dynamics of the system. This dissipativity will give us
estimate of the atmospheric temperature evolution under oceanic feedback.
Dynamical properties that follow  from this dissipativity will be considered in the next section. In particular, we will show that  the coupled  atmosphere-ocean system
has a random attractor, has finite degree of freedom,  and is ergodic
under suitable conditions.\\

We introduce the spaces
\begin{align*}
\tilde H&= L_2(0,1)\times L_2(D)\times\dot{L}_2(D)\\
\tilde V&=W_2^1(0,1)\times W_2^1(D)\times \dot{W}_2^1(D).
\end{align*}
We also choose a subset of  dynamical variables of our system (\ref{eqn1}).
\begin{equation}
\label{eqE}
\tilde v=(\tilde \Theta,T,S),\quad \tilde\Theta=\Theta-z.
\end{equation}
 To calculate the energy inequality
for $\tilde v$,   we apply the chain rule to $\|\tilde v\|_H^2$.
We obtain by Lemma \ref{lcon}
\begin{align}\label{eqE1}
\begin{split}
\frac{d}{dt} &\|\tilde v\|_{\tilde H}^2+2\|\nabla\tilde v\|_{L_2}^2
+2\|b^\frac{1}{2} \tilde\Theta\|_{L_2}^2+2\|\tilde\Theta\|_{L_2}^2\\
=&-2(a,\tilde\Theta)_{L_2}+2(S_a,\tilde\Theta)_{L_2}-
2(b \,S_o,\tilde\Theta)_{L_2}+2(b\,\gamma_{z=1}T,\tilde\Theta)_{L_2}\\
&+2\langle\gamma_{z=1}^\ast S_o,T\rangle+2\langle\gamma_{z=1}^\ast\tilde\Theta,T\rangle+
2\langle\gamma_{z=1}^\ast z(\theta_t\omega),T\rangle
-2\langle\gamma_{z=1}^\ast\gamma_{z=1}T,T\rangle\\
&+2\langle\gamma_{z=1}^\ast F,S\rangle.
\end{split}
\end{align}
Here and in the following we stress  that $0 < b(y) < 1$. The
expression $\nabla \tilde v$ is defined by
$(\partial_y\tilde\Theta,\nabla_{y,z}T,\nabla_{y,z}S)$. We now can
estimate the terms on the right hand side. We have the following
estimate for the second line of (\ref{eqE1}) by the Cauchy-Schwarz
inequality
\begin{align*}
4a^2 &+ \frac{1}{4}\|\tilde\Theta\|_{L_2}^2
+4\|S_a\|_{L_2}^2+\frac{1}{4}\|\tilde\Theta\|_{L_2}^2\\
&+2\|b^\frac{1}{2}
S_o\|_{L_2}^2+\frac{1}{2}\|b^\frac{1}{2}\tilde\Theta\|_{L_2}^2
+\frac{2}{3}\|b^\frac{1}{2}\|_{L_{\infty}}\|\gamma_{z=1}T\|_{L_2}^2+\frac{3}{2}\|b^\frac{1}{2}\tilde\Theta\|_{L_2}^2.
\end{align*}
For the next line we obtain the estimate
\begin{align*}
6\|S_o\|_{L_2}^2+&\frac{1}{6}\|\gamma_{z=1}T\|_{L_2}^2\\
+&\|\tilde\Theta\|_{L_2}^2
+\|\gamma_{z=1}T\|_{L_2}^2+6\|z(\theta_t\omega)\|_{L_2}^2+\frac{1}{6}\|\gamma_{z=1}T\|_{L_2}^2
-2\|\gamma_{z=1}T\|_{L_2}^2.
\end{align*}
Now we can estimate the last line of (\ref{eqE1}). For any $\eps>0$ we can find
an $c_8(\eps)>0$ such that
\[
\eps\|S\|_{W_2^1}^2+c_8(\eps)\|F\|_{L_2}^2 .
\]
Here we also applied the trace theorem
$\|\gamma_{z=1}S\|_{H^\frac{1}{2}}\le c_9\|S\|_{W_2^1}$.
Adding all terms containing $\|\gamma_{z=1}T\|_{L_2}^2$, we see that
 the sum is negative.

Collecting all these estimates, we have
\begin{align*}
\frac{d}{dt}\|\tilde v\|_{\tilde H}^2 &+ 2\|\nabla\tilde
v\|_{L_2}^2 + \frac23(1 -
\|b^\frac12\|_{L_{\infty}})\|\gamma_{z=1}T\|_{L^2}^2 +
\frac12\|\tilde \Theta\|_{L^2}^2
\\
\le 4a^2 &+ 4\|S_a\|_{L_2}^2 + 2\|b^{\frac12}S_o\|_{L_2}^2
+ 6\|S_o\|_{L_2}^2 + 6\|z(\theta_t\omega)\|_{L_2}^2+c_8(\eps)\|F\|_{L_2}^2.
\end{align*}
By using the Poincar\'e inequality for  $S\in \dot{W}_2^1(D)$,
(\ref{tcoest}), and choosing $\eps$ small enough, we conclude that
there is a positive dissipativity constant$\alpha$ such that
 \begin{equation}  \label{energy}
\frac{d}{dt}\|\tilde v\|_{\tilde H}^2 + \alpha(\|\tilde v\|_{\tilde H}^2
+\|\nabla\tilde v\|_{L_2}^2) \le c_{10}+6\|z(\theta_t\omega)\|_{L_2}^2
 \end{equation}
 where $\alpha,\,c_{10}$ is determined by physical  data $a^2, \|S_a\|_{L_2},\,\|S_o\|_{L_2}, \|F\|_{L_2}$
 and $\|b\|_{L^{\infty}}$.\\

By the Gronwall inequality, we finally conclude that

\begin{equation} \label{feedback}
\|\tilde v\|_{\tilde H}^2 \leq
    \|\tilde v(0)\|_{\tilde H}^2e^{-\alpha t}+\frac{c_{10}}{\alpha}
    + 6e^{-\alpha t}\int_0^t \|z(\theta_s\omega)\|_{L_2}^2e^{\alpha s} ds .
\end{equation}

We now show the dissipativity of $\tilde v$ and $v$. Roughly speaking dissipativity means
that all trajectories of the system move to a bounded set in the phase space.
For a random system we have
the following version of dissipativity.
\begin{definition}\label{defA}
A random set $B=\{B(\omega)\}_{\omega\in \Omega}$ consisting of closed bounded sets $B(\omega)$
is called absorbing for a random dynamical system $\phi$ if we have for
any  random set $D=\{D(\omega)\}_{\omega\in\Omega},\,D(\omega)\in H$ bounded,
 such that $t\to \sup_{y\in D(\theta_t\omega)}\|y\|_H$ has a subexponential
growth for $t\to\pm\infty$
\begin{align}\label{eqA0}
\begin{split}
&\phi(t,\omega,D(\omega))\subset B(\theta_t\omega) \quad\text{for }t\ge t_0(D,\omega)\\
&\phi(t,\theta_{-t}\omega,D(\theta_{-t}\omega))\subset B(\omega)\quad\text{for }t\ge t_0(D,\omega).
\end{split}
\end{align}
$B$ is called forward invariant if
\[
\phi(t,\omega,u_0)\in B(\theta_t\omega)\quad \text{if }  u_0\in B(\omega)\quad \text{for }t\ge 0.
\]
\end{definition}
Although $\tilde v$ is not a random dynamical system in the strong sense
we can also show dissipativity in the sense of the above definition.

\begin{lemma}\label{lemA}
Let $\tilde \phi(t,\omega,v_0)\in \tilde H$ for $v_0\in H$ be defined in (\ref{eq19}).
Then the closed ball
$B(0,R_1(\omega))$ with radius
\[
R_1(\omega)= 2\int_{-\infty}^0e^{\alpha\tau}(c_{10}+6\|z(\theta_\tau\omega)\|_{L_2}^2)d\tau
\]
is forward invariant and absorbing.
\end{lemma}

The proof of this lemma follows by integration of (\ref{energy}).\\

It remains to prove the dissipativity of the dynamical system $\phi$.
To this end we obtain from the second equation of
(\ref{eqn1}):
\begin{align*}
\begin{split}
\frac{d}{dt}\|q\|_{L_2}^2+c_{11}^2{\rm Pr}\|q\|_{L_2}^2
&\le {\rm Pr\,{Ra}^2}\|\nabla \tilde v\|_{L_2}^2\\
&\le \frac{{\rm Pr\,{Ra}^2}c_{10}}{\alpha}+6\frac{{\rm Pr\,{Ra}^2}}{\alpha} \|z(\theta_t\omega)\|_{L_2}^2
-\frac{{\rm Pr\,{Ra}^2}}{\alpha}\frac{d}{dt}\|\tilde v\|_{\tilde H}^2
\end{split}
\end{align*}
with some embedding constant $c_{11}$  in the Poincar\'e
inequality $\|q\|_{L_2}\le c_{11}\|\nabla q\|_{L_2}$ for $q\in
\Wsob$. Note that $q$ satisfies homogeneous Dirichlet boundary
conditions. Hence the variation of constants formula allows us to
estimate:
\begin{align}\label{eqA}
\begin{split}
\|q(t,\omega,u_0)\|_{L_2}^2\le&\|v_0\|_H^2e^{-c_{11}^2{\rm Pr}\,t}
+\int_0^t
\frac{({\rm Pr\,{Ra}^2}c_{10}}{\alpha}+6\frac{{\rm Pr\,{Ra}^2}}{\alpha} \|z(\theta_s\omega)\|_{L_2}^2\\
&-\frac{{\rm Pr\,{Ra}^2}}{\alpha}\frac{d}{ds}\|\tilde v(s)\|_{\tilde H}^2)
e^{-c_{11}^2{\rm Pr}\,(t-s)}ds.
\end{split}
\end{align}
Now we apply the integration by parts:
\begin{align*}
-\int_0^t\frac{d}{ds}\|\tilde v(s)\|_{\tilde H}^2e^{c_{11}^2{\rm Pr}\,s}ds=&
\int_0^tc_{11}^2{\rm Pr}\|\tilde v(s)\|_{L_2}^2e^{c_{11}^2{\rm Pr}\,s}ds\\
&-e^{c_{11}^2{\rm Pr}t} \|\tilde v(t)\|_{L_2}^2+
\|\tilde v(0)\|_{L_2}^2.
\end{align*}
Substituting this equation into (\ref{eqA}), we obtain
\begin{align} \label{qqq}
\begin{split}
\|q(t,\omega,u)\|_{L_2}^2 &\le \|v_0\|_H^2e^{-c_{11}^2{\rm Pr}\,t}+
\int_0^t
(\frac{{\rm Pr\,{Ra}^2}}{\alpha}\,c_{10}+6\frac{{\rm Pr\,{Ra}^2}}{\alpha} \|z(\theta_s\omega)\|_{L_2}^2) e^{-c_{11}^2{\rm Pr}(t-s)}ds\\
&+\int_0^t(\frac{c_{11}^2{\rm {Ra}^2}{\rm Pr}^2}{\alpha}\|\tilde v(s)\|_{\tilde H}^2)e^{-c_{11}^2{\rm Pr}(t-s)}ds
+
e^{-{c_{11}^2\rm Pr}\,t}\frac{{\rm Pr\,{Ra}^2}}{\alpha}\|\tilde v(0)\|_{\tilde H}^2.
\end{split}
\end{align}
Note that $\|\tilde v(t)\|_{\tilde H}^2$
is bounded by
\[
\|\tilde v_0\|_{H}^2e^{-\alpha t} + R_1(\theta_t\omega)
\]
which follows from (\ref{energy}).
To construct the radius of the absorbing set we have to replace $\omega$ by $\theta_{-t}\omega$.
Suppose that $t\to\|v_0(\theta_{-t}\omega)\|_{H}^2$ growths not faster than
subexponential.
Then we have that
\[
\lim_{t\to\infty}\|v_0(\omega)\|_{H}^2e^{-\alpha t}=0,\quad
\lim_{t\to\infty}\|v_0(\theta_{-t}\omega)\|_{H}^2e^{-\alpha t}=0.
\]
Hence we can conclude
\[
\lim_{t\to\infty}\int_0^t\frac{c_{11}^2{\rm {Ra}^2}{\rm Pr}^2}{\alpha}\|v_0(\omega)\|_{H}^2e^{-\alpha t}e^{-c_{11}^2{\rm Pr}(t-s)}ds=0
\]
and
\[
\lim_{t\to\infty}\int_0^tc_{11}^2\frac{{\rm {Ra}^2}{\rm Pr}^2}{\alpha}\|v_0(\theta_{-t}\omega)\|_{H}^2e^{-\alpha\,t}e^{-c_{11}^2{\rm Pr}(t-s)}ds=0.
\]
We also note that
\begin{align*}
\lim_{t\to\infty}&\int_0^t(\frac{{\rm Pr\,{Ra}^2}c_{10}}{\alpha}+
6\frac{{\rm Pr\,{Ra}^2}}{\alpha} \|z(\theta_{t-s}\omega)\|_{L_2}^2
+
\frac{c_{11}^2{\rm {Ra}^2}{\rm Pr}^2}{\alpha}R_1(\theta_{s-t}\omega))e^{-{c_{11}^2\rm Pr}(t-s)}ds\\
=\lim_{t\to\infty}&\int_{-t}^0(\frac{{\rm Pr\,{Ra}^2}c_{10}}{\alpha}+
6\frac{{\rm Pr\,{Ra}^2}}{\alpha} \|z(\theta_s\omega)\|_{L_2}^2
+
\frac{c_{11}^2{\rm {Ra}^2}{\rm Pr}^2}{\alpha}R_1(\theta_{s}\omega))e^{{c_{11}^2\rm Pr}s}ds\\
&=:\frac{R_2(\omega)}{2}
<\infty.
\end{align*}

The finiteness of this limit follows because the growth of $t\to R_1(\theta_t\omega)$ is subexponential.
Since all other terms are coupled with exponentially decreasing factors we have found:
\begin{lemma}\label{lEX}
Suppose that the assumptions of Lemma \ref{lemA} are satisfied. Then the random set
$\{B(\omega)\}_{\omega\in\Omega}$
given by closed balls $B(0,R(\omega))$ in $H$
with center zero and radius $R(\omega):=R_1(\omega)+R_2(\omega)$ is an absorbing
and forward invariant set denoted by $B(\omega)$
for the
random dynamical system generated by (\ref{eq18}).
\end{lemma}

For the applications in the next section we need that the elements which are contained in the
absorbing set satisfy a particular regularity. To this end we introduce the function space
\[
{\mathcal{H}}^s:=\{u\in H: \|u\|_s^2:=\|A^\frac{s}{2}u\|_H^2<\infty\}
\]
where $s\in {\mathbb{R}}$. The operator $A^s$ is the $s$-th power of the positive
and symmetric operator $A$.
Note that these spaces are embedded in the Slobodeckij spaces
$H^s,\, s>0$. The norm of these spaces is denoted by $\|\cdot\|_{H^s}$.
This norm can be found in Egorov and Shubin \cite{EgoShu91}, Page 118.
But we do not need this norm explicitly. We only mention that on ${\mathcal{H}}^s$ the norm
$\|\cdot\|_s$
of
$H^s$ is equivalent to the norm of ${\mathcal{H}}^s$ for $0 < s$, see \cite{LioMag68}.\\
The reason to introduce these spaces is that the trace theorem can be formulated with respect to
$H^s$; see Egorov and Shubin \cite{EgoShu91}, Page 120:

\begin{theorem}
Assume that $\alpha>\frac{1}{2}$. Then the  trace mapping between
two Sobolev/Slobodeckij spaces
\[
\gamma_{z=1}: H^\alpha(D) \to  H^{\alpha-\frac{1}{2}}(0,1)
\]
is continuous.
\end{theorem}
This formula generalized (\ref{trace1}) because we can take for $\alpha=1$ and $H^1=W_2^1(D)$
or $V$.
Our goal is it to show that  $v(1,\omega,D)$ is  a bounded set in
${\mathcal{H}}^s$ for some $s>0$.
This property causes the complete continuity of the mapping $v(1,\omega,\cdot)$.
We now derive a differential inequality for $t\|v(t)\|_s^2$.
By the chain rule we have
\[
\frac{d}{dt}(t\|v(t)\|_s^2)=\|v(t)\|_s^2+t\frac{d}{dt}\|v(t)\|_s^2.
\]
Note that for the embedding  constant $c_{12,s}$ between ${\mathcal{H}}^s$ and $V$
\[
\int_0^t\|v\|_s^2ds\le c_{12,s}^2 \int_0^t\|v\|_V^2ds\quad \text{for } s \le 1
\]
such that the left hand  side is bounded if the initial conditions $v_0$
are contained in a bounded set in $H$.
The second term in the above formula can be expressed as followed:
\begin{align*}
t\frac{d}{dt}(A^\frac{s}{2}v,A^\frac{s}{2}v)_H=&2t(\frac{d}{dt}v,A^sv)_H
=-2t(Au,A^sv)_H+2t(F_1(v),A^sv)_H\\
&+2t(F_2(v+Z(\theta_t\omega)),A^sv)_H.
\end{align*}
We have
\[
(Av,A^sv)_H=\|A^{\frac{1}{2}+\frac{s}{2}}v\|_H=\|v\|_{1+s}^2.
\]
If we apply some embedding theorems, see Temam \cite{Tem83} Page 12
we have got for a $c_{13}>0$
\begin{align*}
(F_1(v),\zeta)_H&\le c_{13}\,\|v\|_{m_1+1}\|\psi\|_{m_2+1}\|\zeta\|_{m_3}\,\quad \zeta\in H_{m_3}
\end{align*}
where
$m_1+m_2+m_3\ge 1$ and $0\le m_i<1$. Here we use that $D$ is of dimension 2.
We then have for $m_1=0,\,m_2=s<1$ and $m_3=1-s$
\[
|( F_1(v),A^sv)_H|\le c_{13}\|v\|_V\|\psi\|_{1+s}\|v\|_{1+s}.
\]
$\|\psi(t)\|_{1+s}$ is bounded by $c_1^\prime\|v(t)\|_H$
similar to (\ref{eq6a}) and $\|v(t)\|_{L_\infty(0,T;H)}<\infty$.
To ensure that the norm of $\|\psi\|_{m_2+1}$ is well defined we set $\psi=(0,\psi,0,0)$.
Hence we have for any $\eps>0$ a constant $c_{14}(\eps)$:
\begin{align*}
(F_1(v(t)),A^s v(t))_H&\le c_{14}(\eps) \|v\|_{L_\infty(0,T;H)}^2\|v(t)\|_V^2+\eps \|v(t)\|_{1+s}^2,
\end{align*}
where $\eps$ is chosen sufficiently small.\\
To obtain an estimate for the expression $\langle F_2(v-Z),A^sv\rangle$ we only consider
the expressions
$\langle \gamma_{z=1}^\ast S_o,A^sv\rangle$ and
$\langle \gamma_{z=1}^\ast\gamma_{z=1}T,A^sv\rangle$.
The other term can be treated similarly.
Suppose that $S\in L_2(0,1)$.
We interpret $\gamma^\ast S_o$ as $(0,0,0,\gamma^\ast S_o)$.
Then we have for  any $\eps>0$ a $c_{15}(\eps),\,c_{16}(\eps)>0$ such that for $0<s<\frac{1}{2}$
\begin{align*}
\langle \gamma_{z=1}^\ast S_o,A^sv\rangle&\le \eps\|A^{\frac{1}{2}+\frac{s}{2}}v\|_H^2+
c_{15}(\eps)\|A^{-\frac{1}{2}+\frac{s}{2}}\gamma_{z=1}^\ast S_o\|_H^2\\
&=
\eps\|v\|_{1+s}^2
+c_{15}(\eps)\|\gamma_{z=1}^\ast S_o\|_{-1+s}^2\\
&\le \eps\|v\|_{1+s}^2 + c_{16}(\eps)\|S_o\|_{H^{s-\frac{1}{2}}}^2
\le \eps\|v\|_{1+s}^2+c_{16}(\eps)\|S_o\|_{L_2(0,1)}^2.
\end{align*}
Since by the trace theorem we have
$\gamma_{z=1}^\ast: H^{-\frac{1}{2}+s}(0,1)\to H^{-1+s},\,1-s>\frac{1}{2}$.
For $\alpha<0$ the space $H^\alpha$ denotes the dual space of $H^{-\alpha}$.\\
For $\langle \gamma_{z=1}^\ast\gamma_{z=1}T,A^sv\rangle$ we have for positive constants
$c_{17}$-$c_{19}$, a
sufficiently small
$\eps>0$, $c_{20}(\eps)$ and a sufficiently small $\eps^\prime>0$
\begin{align*}
\langle \gamma_{z=1}^\ast\gamma_{z=1}T,A^sv\rangle=&(\gamma_{z=1}T,\gamma_{z=1}A^sv)_{L_2}\\
&
\le
\|\gamma_{z=1}T\|_{L_2}
\|\gamma_{z=1}A^sv\|_{L_2}\\
&\le
c_{17}\|v\|_V\|\gamma_{z=1}A^sv\|_{L_2}\\
&\le c_{18}\|v\|_V\|A^{\frac{1}{4}+\eps^\prime}A^sv\|_H
\le c_{19}\|v\|_V\|A^{\frac{1}{2}+\frac{s}{2}}v\|_H
\le c_{19}\|v\|_V\|v\|_{1+s}\\
&\le c_{20}(\eps)\|v\|_{V}^2+\eps\|v\|_{1+s}^2.
\end{align*}
For $0<s<\frac{1}{4}$ we have $A^sv\in H^{1-2s}$ for $v\in V$ such that
$\gamma_{z=1}A^s v\in L_2(0,1)$.
Collecting all these estimates we obtain that
$t\|v(t,\omega,v_0)\|_H$ is bounded for $t\le T<\infty$ if $v_0$ is contained in a bounded set.
This allows us to write down the main assertion with respect
to the dissipativity of this section.
\begin{theorem}\label{tA}
For the random dynamical system generated by (\ref{eq18}),
there exists a compact random  set $B=\{B(\omega)\}_{\omega\in \Omega}$
which satisfies Definition \ref{defA}.
\end{theorem}
We define
\begin{equation}\label{eqAB}
B(\omega)=\overline{\phi(1,\theta_{-1}\omega,B(0,R(\theta_{-1}\omega)))}\subset {\mathcal{H}}^s,\quad 0<s<\frac14
\end{equation}
In particular, ${\mathcal{H}}^s$ is compactly embedded in $H$.

\section{Random dynamical behavior}

In this section we will apply the  dissipativity result of the last section to
analyse the dynamical behavior of  the coupled atmosphere-ocean
system (\ref{eqn1}). However, it will
be enough to analyse the transformed  random
dynamical system generated by
(\ref{eq19}). By the transformation (\ref{eq19a}) we can take over
all these qualitative properties to the system (\ref{eq19}).\\

We will consider following dynamical behavior:
random climatic attractors, finite degrees of freedom,
atmospheric temperature evolution under oceanic feedback, and
ergodicity.

We first consider random climatic attractors.
We recall the following  basic concept; see, for instance, Flandoli and Schmalfu{\ss} \cite{FlaSchm95a}.

\begin{definition}
Let $\phi$ be a random dynamical
dynamical system.
A random set $A=\{A(\omega)\}_{\omega\in \Omega}$ consisting of
compact nonempty sets $A(\omega)$ is called random global attractor
if for any random bounded set $D$ we have for the limit in probability
\[
({\mathbb{P}})\lim_{t\to\infty} {\rm dist}_H(\phi(t,\omega,D(\omega)),A(\theta_t\omega)) = 0
\]
and
\[
\phi(t,\omega,A(\omega))=A(\theta_t\omega)
\]
any $t\ge 0$ and $\omega\in\Omega$.
\end{definition}

The essential long-time behavior of a random system
is captured by a random  attractor.
In the last section we showed that the dynamical system $\phi$ generated by
(\ref{eq19}) is dissipative which means that there exists a random set $B$ satisfying
(\ref{eqA0}). In addition,  this set is compact. We now recall and adapt the following
theorem  from \cite{FlaSchm95a}.

\begin{theorem}
Let $\phi$ be a random dynamical
dynamical system on the state space $H$ which is a separable Banach space such that
$x\to\phi(t,\omega,x)$ is continuous.
Suppose that $B$ is a set ensuring the dissipativity given in definition \ref{defA}.
In addition, $B$ has a subexponential growth
(see  Definition \ref{defA})  and is regular (compact).
Then the dynamical system   $\phi$ has
a random attractor.
\end{theorem}

This theorem can be applied to our random dynamical system
 $\phi$ generated by the stochastic differential equation (\ref{eq19}).
 Indeed, all the assumptions
are satisfied. The set $B$ is defined in Theorem \ref{tA}. Its subexponential growth
follows from $B(\omega)\subset B(0,R(\omega))$ where
the radius $R(\omega)$ has been introduced in the last section.
Note that $\phi$ is a {\em continuous} random dynamical system; see Theorem \ref{tEX}. Thus $\phi$ has a random attractor.
By the transformation (\ref{eq19a}), this is also true for the original
coupled atmosphere-ocean system.

\begin{coro}  \label{attractor}
({\bf Random Attractor})
 The coupled atmosphere-ocean system (\ref{eqn1})
has a random attractor.
\end{coro}

\bigskip

 Dissipative systems often have  finite   degrees  of freedom.
This is reflected by the fact that the  Hausdorff-dimension of the attractor is
finite. This fact can be applied to fluid dynamical systems; see for instance
Temam  \cite{Tem97}, Page 403 ff.
A similar theory has been developed for random
dynamical systems; see, for example, \cite{CraFla99, Deb98, Schm97c}.\\

However, we will follow another approach to show that the
random climatic attractor of $(\ref{eqn1})$ has
only finitely many degrees of freedom. Namely we will use the technique of determining
functionals. This technique has been introduced for deterministic systems by Foias and Prodi \cite{FoiPro67}, and
Ladyzhenskya \cite{La2}.  See   \cite{JohTit93, C3} for more recent work.
 Roughly speaking,
a set of determining functionals is a set of functionals (for instance, Fourier modes)
such that if it is known that a dynamical system has an asymptotic stable behavior
only with respect to these finitely many modes, then the complete system has an asymptotic
stable behavior.
For random dynamical systems we can investigate
{\em determining functionals in probability} as in  Chueshov et al. \cite{ChuDuaSchm01a}.

    \begin{definition}
   We call a set ${\mathcal{L}}=\{l_j,\;j=1,\cdots, N\}$ of linear continuous
   and linearly independent functionals on a space $X$ continuously embedded in $H$
   (for instance $X={\mathcal{H}}^s$ or $V$)
   asymptotically determining in probability
   if
   \[
   ({\mathbb{P}})\lim_{t\to\infty}\int_t^{t+1}
   \max_j|l_j(\phi(\tau,\omega,v_1)-\phi(\tau,\omega,v_2))|^2d\tau= 0
   \]
   for two initial conditions $v_1,\,v_2\in H$ implies
   \[
   ({\mathbb{P}})\lim_{t\to\infty}\|\phi(t,\omega,v_1)-\phi(t,\omega,v_2)\|_H=
   0.
   \]
   \end{definition}

Often  the  elements of set ${\mathcal L}$  can be chosen
as the projections with respect to
the Fourier expansion of the solution. Since our
domain $D$ is an rectangle, we can calculate the
Fourier expansion more explicitly.
\\
In the following we need an additive embedding inequality     based on
qualitative difference of the spaces $H$ and $X$ for some set of functionals
${\mathcal{L}}$
\begin{equation}\label{eqe8}
\|u\|_H\le C_{{\mathcal{L}}}\max_{l_i\in{\mathcal{L}}}|l_i(u)|+\eps_{{\mathcal{L}}}\| u\|_X,\quad C_{\mathcal{L}}>0
\end{equation}
The constant $\eps_{{\mathcal{L}}}>0$ describes a fundamental difference of the spaces $X$ and $H$.
How (\ref{eqe8}) works is described by an motivating example
in Chueshov et al. \cite{ChuDuaSchm01a}.

 We recall and adapt a result from \cite{ChuDuaSchm01a}.

\begin{theorem}\label{t1}
Let ${\mathcal L} =\{ l_j : j= 1,...,N\}$ be  a set of linear
continuous and linearly independent functionals on $X$. We assume
that we have an absorbing and forward invariant set $B$ in $X$
such that for $\sup_{v\in B(\omega)}\|v\|_{X}^2$ the expression
$t\to\sup_{v\in B(\theta_t\omega)}\|v\|_{X}^2$ is locally integrable
and subexponentially growing. Suppose there exist a constant
$c_{21}>0$ and a measurable function $l\ge0$ such that for
$v_1,\,v_2\in V$ we have for
$F(\omega,v)=F_1(v)+F_2(v+Z(\omega))$
\begin{align*}
\begin{split}
\langle - A(v_1-v_2)+
F(\omega,v_1)-F(\omega,v_2),v_1-v_2\rangle
\\
\le
-c_{21}\|v_1-v_2\|_V^2+ l(v_1,v_2,\omega)\|v_1-v_2\|_H^2.
\end{split}
\end{align*}
Assume that
\begin{equation*}
\frac{1}{m}{\mathbb{E}}\left\{\sup_{v_1,v_2\in B(\omega)}\int_0^m
l(\phi(t,\omega,v_1),\phi(t,\omega,v_2),\theta_t\omega)dt\right\}<
c_{21}\eps_{{\mathcal L}}^{-2}
\end{equation*}
for some $m>0$. Then ${\mathcal L}$ is a set of asymptotically
determining functionals in probability for random dynamical system
$\phi$.
\end{theorem}

We set $X={\mathcal{H}}^s,  s\in (0,\frac14)$.
In the last section we have shown that the set $B$,
consisting of {\em bounded} sets,
is forward invariant. The function
$l$ appearing in the formulation of Theorem \ref{t1} expresses the local Lipschitz
continuity of our nonlinear operator $F(v)$. The essential part that determines $F$ is
the Jacobian operator $J$.  A method  on how to estimate the function $l$
is  in \cite{ChuDuaSchm01a}. We also note the local Lipschitz constant can be
estimated in terms of the ${\mathcal{H}}^s$-norm.
To estimate $F_2$ we have to apply the techniques
introduced in Section 4. Note that $F_2$ is linear.
\\
We  can apply  Theorem \ref{t1} to   the random dynamical system
generated by (\ref{eq19}) and get the following result.

\begin{theorem}
The random dynamical system generated by (\ref{eq19}) has finitely many
degrees of freedom. More precisely, there exists a set of
      linearly   independent continuous
functionals (on ${\mathcal{H}}^s$) which is
asymptotically determining in probability.
\end{theorem}

Because of  the transformation (\ref{eq19a}), this result is also true
for the coupled atmosphere-ocean system:

\begin{coro}   \label{finite}
({\bf Finite Degrees of Freedom})
The   coupled atmosphere-ocean system (\ref{eqn1}) has finitely many
degrees of freedom, in the sense of having a finite set of
 linearly independent   continuous functionals  which is
asymptotically determining in probability.
\end{coro}

The number of the elements in ${\mathcal{L}}$ can be estimated explicitly as  in \cite{ChuDuaSchm01a}.\\

\bigskip

Now we   consider random fixed point and ergodicity.
We can do a small modification of (\ref{eqn1}). This modification is given when we
replace $\Delta T$ by $\nu\Delta T$ and $\Delta S$ by $\nu\Delta S$ where $\nu>0$ is viscosity.
Under particular assumptions about physical data in (\ref{eqn1}) ,
we can show that the
behavior of our dynamical system is laminar.
For a stochastic system,  this
means that after a relatively short time,
all trajectories starting from different initial states show almost the same
dynamical behavior.
This can be seen easily if $S_o,\,a,\,S_a,\,F$ are zero, there is no
noise and $\nu$ is large.
We will show that a laminar behavior also appears when  $S_o,\,a,\,S_a,\,F$  are small in some sense.\\
Mathematically speaking,  laminar behavior means that a random dynamical system has a unique   exponentially
attracting random fixed point.

\begin{definition}
A random variable $v^\ast:\Omega\to H$ is defined to be a random fixed point for a random dynamical system
if
\[
\phi(t,\omega,v^\ast(\omega))=v^\ast(\theta_t\omega)
\]
for $t\ge 0$ and $\omega\in\Omega$. A random fixed point
$v^\ast$ is called exponentially attracting
if
\[
\lim_{t\to\infty}\|\phi(t,\omega,x)-v^\ast(\theta_t\omega)\|_H=0
\]
for any $x\in H$ and $\omega\in \Omega$.
\end{definition}
Sufficient conditions for the existence of random fixed points
are given in Schmalfu{\ss}
\cite{Schm97a}. We here formulate a simpler version of this
theorem and it is appropriate  for our system here.

\begin{theorem} ({\bf Random Fixed Point Theorem})
Let $\phi$ be a random dynamical system and suppose that $B$ is a forward invariant complete set.
In addition, $B$ has a subexponential growth, see Definition \ref{defA}.
Suppose that the following contraction conditions holds:
\begin{equation}\label{eqKX}
\sup_{v_1\not= v_2\in B(\omega)}\frac{\|\phi(1,\omega,v_1)-\phi(1,\omega,v_2)\|_H}{\|v_1-v_2\|_H}\le k(\omega)
\end{equation}
where the expectation of $\log k$ denoted by ${\mathbb{E}}\log k<0$.
Then $\phi$ has a unique random fixed point in $B$ which is exponentially attracting.
\end{theorem}

This theorem can be considered as a random version of the Banach fixed point theorem.
The contraction condition is formulated in the mean for the right hand side of  (\ref{eqKX}).\\

\begin{theorem} \label{fpt}
Assume  that the physical data $|a|,\,\|S_a\|_{L_2},\,\|S_o\|_{L_2},\,\|F\|_{L_2}$ and
the trace of the covariance for the noise ${\rm tr}_HQ$
are sufficiently small,  and that the viscosity $\nu$ is sufficiently large.
Then the random dynamical system generated by (\ref{eq19})
has a unique random fixed point in $B$.
\end{theorem}
Here we only give a short sketch of the proof.
Let us suppose for a while that $B$ is given by the ball $B(0,R)$ introduced in Lemma \ref{lEX}.
Suppose that the data in the assumption of the lemma are small and $\nu$ is large.
Then it follows that ${\mathbb{E}}R$ is also small.
To calculate the contraction condition we have to calculate
 $\|\phi(1,\omega,v_1(\omega))-\phi(1,\omega,v_2(\omega))\|_H^2$
for arbitrary random variables  $v_1,\,v_2\in B$.
By Lemma \ref{lcon} we have that
\[
\langle J(q_1,\psi_1)-J(q_1,\psi_1),q_1-q_2\rangle\le
c_{22}\|q_1-q_2\|_{W_2^1}^2+c_{23}\|q_1\|_{\Wsob}^2\|q_1-q_2\|_{L_2}^2
\]
where the constant $c_{23}$ can be chosen sufficiently small if $\nu$ is large. On account
of the fact that the other expressions allow similar estimates and that $F_2$ is linear
we obtain:
\begin{align*}
\frac{d}{dt}&\|\phi(t,\omega,v_1(\omega))-\phi(t,\omega,v_2(\omega))\|_H^2\\
&\le
(-\alpha^\prime+c_{23}\|\phi(t,\omega,v_2(\omega))\|_V^2)\|\phi(t,\omega,v_1(\omega))-\phi(t,\omega,v_2(\omega))\|_H^2
\end{align*}
for some positive $\alpha^\prime $ depending on $\nu$.
From this inequality and the Gronwall lemma it follows that the contraction condition (\ref{eqKX}) is satisfied if
\[
{\mathbb{E}}\sup_{u_2\in
B(\omega)}c_{23}\int_0^1\|\phi(t,\omega,v_2)\|_V^2dt<\alpha^\prime.
\]
But by the energy inequality this property is satisfied if the ${\mathbb{E}}R$ and
${\mathbb{E}}\|z\|_V^2$ is sufficiently small which follows from the assumptions.
\\

Let now $B$ be the random set defined in (\ref{eqAB}).
Since the set $B$ introduced in (\ref{eqAB}) is absorbing any state the fixed point
$v^\ast$ is contained in this $B$. In addition $v^\ast$ attracts {\em any} state
from $H$ and not only states from $B$.

\begin{coro}\label{corrf}
({\bf Unique Random Fixed Point}) Assume that the  physical data
$|a|,\,\|S_a\|_{L_2},\,\|S_o\|_{L_2},\,\|F\|_{L_2}$ and the trace
of the covariance for the noise ${\rm tr}_HQ$ are sufficiently
small, and that the viscosity $\nu$ is sufficiently large. Then,
through the transformation (\ref{eq19a}), the original coupled
atmosphere-ocean system (\ref{eqn1}) has a unique exponentially
attracting random fixed point
$u^\ast(\omega)=v^\ast(\omega)+Z(\omega)$, where
$u=(\Theta, q, T,S)$.
\end{coro}

The uniqueness of  this random fixed point implies {\em ergodicity}.
We will comment on this issue at the end of this section.

\bigskip

By the well-posedness  Theorem \ref{tEX}, we know that the stochastic
differential equation (\ref{eq19}) for the coupled atmosphere-ocean system  has  a unique solution.  The solution
is a Markov process.  We can define the associated
Markov operators $\mathcal T(t)$ for $t \geq 0$,
as discussed in \cite{Schm91a, Schm89}.
Moreover, $\{{\mathcal T}(t)\}_{t\geq 0}$ forms a semigroup.

Let $M^2$ be the set of probability distributions $\mu$ with finite energy,
i.e.,
\[
\int_H\|u\|_H^2d\mu(u)<\infty .
\]
Then the distribution of the solution $u(t)$ (at time $t$) of the stochastic
differential equation   (\ref{eq19}) is given
by
\[
{\mathcal T}(t)\mu_0,\
\]
where the distribution $\mu_0$ of the initial data is contained in $M^2$.

We note that
the expectation of the solution $\|u(t)\|_H^2$ can be expressed in terms of this
distribution ${\mathcal T}(t)\mu_0$:
\[
{\mathbb{E}}\|u(t)\|_H^2=\int_H \|u\|_H^2d{\mathcal T}(t)\mu_0.
\]
We can derive the following energy inequality in the mean, using our earlier estimates
in (\ref{feedback}) and (\ref{qqq}):

\begin{theorem}\label{tMar}
The dynamical quantity  $u=(\Theta, q, T,S)$ of the coupled
atmospheric-ocean system satisfy the estimate
\[
{\mathbb{E}}\|u(t)\|_H^2+\alpha{\mathbb{E}}\int_0^t\|u(\tau)\|_V^2d\tau
\le {\mathbb{E}}\|u_0\|_H^2+t\,c_{24}+t\,{\rm tr}_{L_2}Q,
\]
where the positive constants $c_{24}$ and $\alpha$
depend    on   physical data  $F(y)$,
$a$,  $S_a(y)$, $ S_o(y)$,  $Pr$ and  $Ra$.
\end{theorem}

By the Gronwall inequality, we further obtain the following result
about the asymptotic mean-square estimate for the coupled
atmosphere-ocean system.

\begin{coro} \label{feedbacktheorem}
 ({\bf Atmospheric Temperature Evolution under Oceanic Feedback})
For the expectation of the dynamical quantity  $u=(\Theta, q, T,S)$ of the coupled
atmospheric-ocean system, we have the asymptotic estimate
 \[
\limsup_{t\to\infty}
{\mathbb{E}}\|u(t)\|_H^2=\limsup_{t\to\infty}\int_H\|u\|_H^2d{\mathcal  T}(t)\mu_0\le
\frac{c_{24}+{\rm tr}_{L_2}Q}{c_{25}}
\]
if the initial distribution $\mu_0$
of the random initial condition $u_0(\omega)$ is contained in
$M^2$.  Here   $c_{25}>0$ also depends on physical data.
In particular, we have asymptotic mean-square estimate for the
atmospheric temperature evolution under oceanic feedback
\begin{equation} \label{temperature}
 \limsup_{t\to\infty}
{\mathbb{E}}\| \Theta \|_H^2
 \le
\frac{c_{24}+{\rm tr}_{L_2}Q}{c_{25}}.
\end{equation}
Thus the atmospheric temperature   $\Theta(y,t)$,
as modeled by the coupled atmosphere-ocean system
 (\ref{eqn1}),
 is bounded  asymptotically   in mean-square norm in terms of
physical quantities such as the freshwater flux $F(y)$,
the trace of the covariance operator of the external noise,
  the earth's longwave radiative cooling coefficient
$a$,
and  the empirical functions $S_a(y)$ and $ S_o(y)$
representing the latitudinal dependence of the shortwave solar radiation,
as well as the  Prandtl  number ${\rm Pr}$ and  the Rayleigh number ${\rm Ra}$ for oceanic fluids.
\end{coro}

By the estimates of Theorem \ref{tMar},
 we are able to use the well known
Krylov-Bogolyubov procedure to conclude the existence
of invariant measures of the Markov semigroup.
\begin{coro}
The semigroup of Markov operators $\{{\mathcal T}(t)\}_{t\geq 0}$ possesses an
invariant distribution $\mu_i$
in $M^2$:
\[
{\mathcal T}(t)\mu_i=\mu_i\quad \text{for }t\ge 0.
\]
\end{coro}

In fact,  the limit points of
\[
\left\{
\frac1t\int_0^t{\mathcal  T}(\tau)\mu_0d\tau\right\}_{t\geq 0}
\]
for $t\to\infty$ are invariant distributions.
The existence of such limit points follows from the estimate
in Theorem \ref{tMar}.

In some situations, the invariant measure may be unique.
For example,
the  unique random fixed point  in  Corollary \ref{corrf} is
 defined by a random variable
$u^\ast(\omega)=v^\ast(\omega)+z(\omega)$.
This random variable corresponds to a unique
invariant measure of the Markov semigroup.   More specifically,
 this unique  invariant measure  is the expectation of the Dirac measure
with the random variable
as the random mass point
\[
\mu_i={\mathbb{E}}\delta_{u^\ast(\omega)}.
\]
Because the uniqueness of   invariant measure
implies   ergodicity \cite{DaPZab96}, we conclude that the
   coupled atmosphere-ocean model (\ref{eqn1})
is ergodic under the   suitable conditions  in  Corollary \ref{corrf} for physical data and random noise.  We reformulate Corollary \ref{corrf}
as the following  ergodicity principle.

\begin{theorem}  \label{ergodic}
({\bf Ergodicity})
Assume that the physical data $|a|,\,\|S_a\|_{L_2},\,\|S_o\|_{L_2},$ $\|F\|_{L_2}$ and
the trace of the covariance for the noise ${\rm tr}_HQ$
are sufficiently small,  and that the viscosity $\nu$ is sufficiently large.
Then the  coupled atmosphere-ocean system  (\ref{eqn1})
is ergodic, namely,   for any observable of
the    coupled  atmosphere-ocean flows, its time average
approximates the statistical ensemble average, as long as  the time
interval is sufficiently long.
\end{theorem}

\section{Summary}

We have investigated the dynamical behavior of a
 coupled atmosphere-ocean model.
First,  we have shown that the asymptotic dynamics
of the coupled atmosphere-ocean model  is described by a
random climatic attractor (Corollary \ref{attractor}).
 Second,  we have estimated the
atmospheric temperature  evolution  under oceanic feedback,
 in terms of the freshwater flux, heat flux
and the   external fluctuation at the air-sea interface, as well as
the earth's longwave radiation coefficient
and the shortwave solar radiation profile (Corollary \ref{feedbacktheorem}).
Third, we have demonstrated that this system has finite
degree of freedom by presenting a finite set of determining
functionals in probability (Corollary \ref{finite}).
Finally, we have proved that the coupled atmosphere-ocean model
is ergodic under suitable conditions for physical parameters and randomness,
and thus for any observable of
the    coupled  atmosphere-ocean flows, its time average
approximates the statistical ensemble average, as long as  the time
interval is sufficiently long (Theorem \ref{ergodic}).

\bigskip

{\bf Acknowledgement.}
This work was partly supported by the NSF Grants DMS-9973204
and DMS-0139073,      the
Grant 10001018 of the NNSF of China, and the Grant BK2001108 of the  NSF
of Jiangsu Province, as well as the
Scientific Research Foundation for Returned Overseas Chinese Scholars of Jiangsu Education Commission.
B. Schmalfuss and H. Gao  would like to thank
Illinois Institute of Technology, Chicago,  and H. Gao  would like to thank
the Institute for Mathematics and Its Applications, Minneapolis,
 for their hospitality.
This research was supported in part
       by the Institute for Mathematics and its Applications with funds
       provided by the National Science Foundation.


\end{document}